\documentclass[10pt,leqno]{amsart}
\usepackage{graphicx}
\usepackage{indentfirst,csquotes}

\usepackage{amssymb,amsthm,amsmath}
\usepackage{xcolor,paralist,hyperref,titlesec,fancyhdr,etoolbox}
\newtheorem{theorem}{Theorem}[]

\newtheorem{example}[theorem]{Example}
\newtheorem{lemma}[theorem]{Lemma}

\titleformat{name=\section}{}{\thetitle.}{0.8em}{\centering\scshape}

\hypersetup{ colorlinks=true, linkcolor=black, filecolor=black, urlcolor=black }

\usepackage{lipsum}

\begin{document}
	\title{Some Parseval-Goldstein Type Theorems For Generalized Integral Transforms} 
	\author[Albayrak]{Durmuş ALBAYRAK}
	\date{27.02.2023}
	\address{Marmara University, Department of Mathematics, 34722, İstanbul-Türkiye}
	\email{durmus.albayrak@marmara.edu.tr}
	\maketitle
	
	\let\thefootnote\relax
	\footnotetext{\textit{2020 Mathematics Subject Classification.} 44A05, 44A10, 44A15, 44A20.}
	\footnotetext{\textit{Key words and phrases.} Generalized Laplace transform; Generalized Stieltjes transform; Laplace transform; Parseval-Goldstein theorem.} 
	
	\begin{abstract}
	In this work, we establish some Parseval-Goldstein type identities and relations that include various new generalized integral transforms such as $\mathcal{L}_{\alpha,\mu}$-transform and generalized Stieltjes transform. In addition, we evaluated improper integrals of some fundamental and special functions using our results.
	\end{abstract} 
	
	\bigskip
	
	\noindent 
	
	\section{Introduction, Definitions and Preliminaries}
	The theory of special functions and integral transforms constitute an important part of research subjects in mathematics, physics and engineering. Generally, an integral transform is defined by
	\begin{align}
		T(y)=\mathcal{T}\left\{f(t);y\right\}=\int_a^b K(y,t)f(t)dt\label{T}
	\end{align} 
	where the function $f(t)$ defined in $a\leq t\leq b$, $K(y,t)$ is called the kernel of transform, and $y$ is called the transform variable \cite{DB}. In the literature, some famous integral transforms are Laplace, Fourier and Stieltjes transforms. Many researchers have defined new integral transforms in the form of (\ref{T}) by choosing different kernels and boundaries. In particular, the kernels of the transforms can be selected from special functions as well as elementary functions. The reader may refer to \cite{DB}.
	
	The Stieltjes transform of a function is obtained by applying the Laplace transform of the function twice. These kinds of relations, where consecutive integral transforms are applied, are referred to as Parseval-Goldstein type relations or theorems. Thus, the image of a function under an unknown new integral transform can be obtained through the successive applications of known integral transforms. As a result, these relations shed light on the calculation of many generalized integrals that have not yet been evaluated.
	
	In 1989, Yürekli \cite{Y1} proved a Parseval-Goldstien type theorem which gives the relationship between Laplace and Stieltjes transforms and many results arising from this theorem. In 1992, he did a similar study for the generalized Stieltjes transform \cite{Y2}. Later, many authors examined similar relationships between different integral transforms based on Parseval-Goldstein type theorems.\cite{A2,A3,KKU,KAU,Y1,Y2}. 
	
	Albayrak \cite{A1} considers a different generalization of Laplace transform over the set of functions
	\begin{equation}
		A=\left\{f(t)|\exists K,M,a\in\mathbb{R}, 
		\left\vert t^{\alpha-\mu}f\left(  t\right)  \right\vert \leq Ke^{at^{\mu}}\text{ for all }t\geq M, K>0\right\},
	\end{equation}
	which is defined by
	\begin{equation}
		F\left(  y\right)  =\mathcal{L}_{\alpha,\mu}\left\{  f\left(  t\right)
		;y\right\}  =\int_{0}^{\infty}t^{\alpha-1}e^{-y^{\mu}t^{\mu}}f\left(
		t\right)  dt,\label{LAM}%
	\end{equation}
	and the inverse of $\mathcal{L}_{\alpha,\mu}%
	-$transform is defined by
	\begin{equation*}
		f(t)=\mathcal{L}_{\alpha,\mu}^{-1}\left\{  F\left(  y\right)  ;t\right\}
		=\dfrac{\mu t^{\mu-\alpha}}{2\pi i}\int_{C}e^{yt^{\mu}}F\left(  y^{1/\mu
		}\right)  dy, \label{ILAM1}%
	\end{equation*}
	where $\alpha,y\in\mathbb{C}$, $\mu\in\mathbb{R}$,  $\operatorname{Re}\alpha> \mu>0$, $\operatorname{Re}y>0$. A generalization of the harmonic oscillator in non-resisting and resisting medium problems, initial-boundary value problems and integral equations are solved via this integral transform. Furthermore, the alternative solution of well-known series entitled as Basel problem is obtained in a similar way. The reader may refer to \cite{A1} for detailed information.
	
	In this study, Parseval-Goldstein type theorem involving $\mathcal{L}_{\alpha,\mu}-$transform will be proved. Later, some generalized integrals will be evaulated as applications of these theorems.
	
	With special choices of $\alpha$ and $\mu$, $\mathcal{L}_{\alpha,\mu}-$transform can be reduced to some classical
	integral transforms, such as $\mathcal{L}_{1,1}\left\{  f\left(  t\right)  ;y\right\}
	=\mathcal{L}\left\{  f\left(  t\right)  ;y\right\}  $ Laplace transform \cite{DB},
	$\mathcal{L}_{2,2}\left\{  f\left(  t\right)  ;y\right\}  =\mathcal{L}%
	_{2}\left\{  f\left(  t\right)  ;y\right\}  $ $\mathcal{L}_{2}-$transform which was introduced by Yürekli and Sadek \cite{YS}, $\mathcal{L}_{\alpha,1}\left\{  f\left(  t\right)  ;y\right\}
	=\mathcal{L}_{\alpha}\left\{  f\left(  t\right)  ;y\right\}$ another generalized Laplace transform which is defined by Karataş et al \cite{KKU,KAU} and
	$\mathcal{L}_{\mu \omega,\mu}\left\{  f\left(  t\right)  ;y\right\}
	=\frac{1}{\omega y^{\mu \omega-1}}\mathfrak{B}_{\omega, \mu}\left\{  f\left(  t\right)  ;y\right\}$ Borel-Džrbashjan transform \cite{BD1,BD2}. If we make a change of variable $t=u^{\frac{1}{\mu}}$ in the right-hand side of (\ref{LAM}), we get
	the following relationship between the Laplace transform and the $\mathcal{L}_{\alpha,\mu}-$transform
	\begin{equation}
		\mathcal{L}_{\alpha,\mu}\left\{  f\left(  t\right)  ;y\right\}  =\dfrac{1}%
		{\mu}\mathcal{L}\left\{  t^{\frac{\alpha}{\mu}-1}f\left(  t^{\frac{1}{\mu}}\right)  ;y^{\mu
		}\right\}  . \label{LAML}%
	\end{equation}
	In the literature, some generalizations of the Stieltjes transform have been examined by many authors and their applications have been included. We will also describe a new generalized Stieltjes transform obtained by applying $\mathcal{L}_{\alpha,\mu}$-transform sequentially. In addition, under appropriate conditions of convergence, we will introduce some new generalized integral transforms with the help of $\mathcal{L}_{\alpha,\mu}$ or some integral transforms.
	
	The generalized Stieltjes-type transform of $f\left(  x\right)  ,$ is defined
	by 
	\begin{equation}
		S_{\alpha,\mu,\rho}(y)=\mathcal{S}_{\alpha,\mu,\rho}\left\{  f\left(  t\right)  ;y\right\}  =\int
		_{0}^{\infty}\dfrac{t^{\alpha-1}}{\left(  y^{\mu}+t^{\mu}\right)  ^{\rho}%
		}f\left(  t\right)  dt \label{SAMR}%
	\end{equation}
	where $\operatorname{Re}\alpha>0,$ $\operatorname{Re}\mu>0,$
	$\operatorname{Re}\rho>0$ and the inverse of generalized Stieltjes-type transform is defined by
	\begin{align*}
		f(t)=\mathcal{S}^{-1}_{\alpha,\mu,\rho}\left\{  S_{\alpha,\mu,\rho}(y)  ;t\right\}=\dfrac{\mu t^{\mu-\alpha}(\rho-1)}{2\pi i}\int_{C}\left(t^\mu+y\right)^{\rho-2}	S_{\alpha,\mu,\rho}\left(  y^{1/\mu
		}\right)  dy,
	\end{align*}
	where $\operatorname{Re}\alpha>0,$ $\operatorname{Re}\mu>0,$
	$\operatorname{Re}\rho>1$. 
	With special choices of $\alpha,\mu$ and $\rho$, $\mathcal{S}_{\alpha,\mu,\rho}-$transform can be reduced to some classical
	integral transforms, such as $\mathcal{S}_{1,1,1}\left\{  f\left(  t\right)
	;y\right\}  =\mathcal{S}\left\{  f\left(  t\right)  ;y\right\}  $ Stieltjes
	transform \cite{DB,E2}, $\mathcal{S}_{2,2,1}\left\{  f\left(  t\right)  ;y\right\}
	=\mathcal{P}\left\{  f\left(  t\right)  ;y\right\}  $ Widder-Potential
	transform \cite{W}, $\mathcal{S}_{1,2,1/2}\left\{  f\left(  t\right)  ;y\right\}
	=\mathcal{G}\left\{  f\left(  t\right)  ;y\right\}$, Glasser transform \cite{G},
	$\mathcal{S}_{1,1,\rho}\left\{  f\left(  t\right)  ;y\right\}  =\mathcal{S}%
	_{\rho}\left\{  f\left(  t\right)  ;y\right\}$ generalized Stieltjes
	transform \cite{E2}, $\mathcal{S}_{2,2,\rho}\left\{  f\left(  t\right)  ;y\right\}
	=\mathcal{P}_{\rho,2}\left\{  f\left(  t\right)  ;y\right\}$ generalized
	Widder-Potential transform\cite{NVYO}.
	
	If we make a change of variable $t=u^{\frac{1}{\mu}}$ in the right-hand side of (\ref{SAMR}), we
	have the following relationship between the generalized Stieltjes transform
	and the generalized Stieltjes-type transform%
	\begin{equation}
		\mathcal{S}_{\alpha,\mu,\rho}\left\{  f\left(  t\right)  ;y\right\}
		=\dfrac{1}{\mu}\mathcal{S}_{\rho}\left\{  t^{\frac{\alpha}{\mu}-1}f\left(  t^{\frac{1}{\mu}}\right)  ;y^{\mu
		}\right\}. \label{SAMRSR}%
	\end{equation}
	Beside the generalized Stieltjes type integral transform, some generalized integral transforms that will shed light on the study will be defined as follows under appropriate convergence conditions.
	
	First, let's give definitions of integral transforms that we want to generalize. Fourier sine and Fourier cosine integral transforms \cite{DB}, respectively, are defined by
	\begin{align}
		\mathcal{F}_{s}\left\{  f\left(  t\right)  ;y\right\}   &
		=\int_{0}^{\infty}\sin\left(  yt\right)  f\left(
		t\right)  dt,\label{FS}\\
		\mathcal{F}_{c}\left\{  f\left(  t\right)  ;y\right\}   &
		=\int_{0}^{\infty}\cos\left(  yt\right)  f\left(
		t\right)  dt.\label{FC}
	\end{align}
	
	Now, we will define a generalized form of these integral transforms under appropriate convergence conditions. 
	
	Generalized Fourier sine and cosine integral transforms are defined as follow
	\begin{align}
		{F}_{s,\alpha,\mu}(y)=\mathcal{F}_{s,\alpha,\mu}\left\{  f\left(  t\right)  ;y\right\}   &
		=\int_{0}^{\infty}t^{\alpha-1}\sin\left(  y^{\mu}t^{\mu}\right)  f\left(
		t\right)  dt,\label{FSAM}\\
		{F}_{c,\alpha,\mu}(y)=\mathcal{F}_{c,\alpha,\mu}\left\{  f\left(  t\right)  ;y\right\}   &
		=\int_{0}^{\infty}t^{\alpha-1}\cos\left(  y^{\mu}t^{\mu}\right)  f\left(
		t\right)  dt,\label{FCAM}
	\end{align}
	where $y^\mu>0$, $f(t)$ is piecewise continuous and $t^{\alpha-1}f(t)$ is absolutely integrable over $\left[0,\infty\right)$. The inverse of $\mathcal{F}_{s,\alpha,\mu}%
	-$transform and $\mathcal{F}_{c,\alpha,\mu}%
	-$transform are defined by
	\begin{align*}
		f(t)=\mathcal{F}^{-1}_{s,\alpha,\mu}\left\{  F_{s,\alpha,\mu}(y)  ;t\right\}   &
		=\frac{2\mu  t^{\mu-\alpha}}{\pi}\int_{0}^{\infty}\sin\left( t^{\mu} y\right)  {F}_{s,\alpha,\mu}(y^{1/\mu})  dy,\\
		f(t)=\mathcal{F}^{-1}_{c,\alpha,\mu}\left\{  F_{c,\alpha,\mu}(y^{1/\mu})  ;t\right\} &
		=\frac{2\mu  t^{\mu-\alpha}}{\pi}\int_{0}^{\infty}\cos\left( t^{\mu} y\right)  {F}_{c,\alpha,\mu}(y^{1/\mu})  dy,
	\end{align*}
	where  ${F}_{s,\alpha,\mu}(y^{1/\mu})$ and ${F}_{c,\alpha,\mu}(y^{1/\mu})$ are piecewise continuous and absolutely integrable over $\left[0,\infty\right)$. Now, let's give special functions that will be used throughout the study \cite{OSM}.
	
	The Gamma function is defined by
	$$
	\Gamma(z)=\int_0^{\infty} t^{z-1} e^{-t} d t, \quad \operatorname{Re}(z)>0 .
	$$
	Basic properties of Gamma function are given in \cite{OSM}.
	Pochammer symbol is defined by the following relation,
	$$
	(\alpha)_n=\left\{\begin{array}{cl}
		\alpha(\alpha+1) \ldots(\alpha+n-1), & n=1,2,3 \ldots \\
		1, & n=0
	\end{array}\right.
	$$
	where $\alpha \in \mathbb{R}$.
	The relationship between the Pochammer symbol and the gamma function is given by
	$$
	(\alpha)_n=\frac{\Gamma(\alpha+n)}{\Gamma(\alpha)}, \quad \alpha \neq 0,1,2, \ldots
	$$
	The generalized hypergeometric series is defined as
	\begin{align}
		{ }_r F_s\left[\begin{array}{c}
			\alpha_1, \alpha_2, \ldots, \alpha_r \\
			\beta_1, \beta_2, \ldots, \beta_s
		\end{array} \mid z\right]=\sum_{n=0}^{\infty} \frac{\left(\alpha_1\right)_n\left(\alpha_2\right)_n \ldots\left(\alpha_r\right)_n}{\left(\beta_1\right)_n\left(\beta_2\right)_n \ldots\left(\beta_s\right)_n} \frac{z^n}{n !},\label{HS}
	\end{align}
	where $r, s \in \mathbb{Z}^{+} \cup\{0\}$ and $\alpha_i, \beta_j \neq 0,-1,-2, \ldots(1 \leq i \leq r, 1 \leq j \leq s)$.  The reader may refer to \cite{OSM} for detailed information about the convergence conditions of this series. The Laplace transform of a generalized hypergeometric function ${ }_r F_s$ in \cite[p.219,Entry(17)]{E1} as follows:
	\begin{align}
		\int_0^{\infty} e^{-yt} t^{v-1}{ }_r F_s\left[\begin{array}{c}
			\alpha_1, \cdots, \alpha_r \\
			\beta_1, \cdots, \beta_s
		\end{array} \mid a t\right] d t=\frac{\Gamma(v)}{y^{v}}{ }_{r+1} F_s\left[\begin{array}{c}
			v, \alpha_1, \cdots, \alpha_r \\
			\beta_1, \cdots, \beta_s
		\end{array} \mid \frac{a}{y}\right]\label{FI}
	\end{align}
	provided if $r<s, \operatorname{Re}(v)>0, \operatorname{Re}(y)>0$ and $a$ is arbitrary or if $r=s>0, \operatorname{Re}(v)>0$ and $\operatorname{Re}(y)>\operatorname{Re}(a)$.
	
	The confluent hypergeometric function is defined \cite{OSM} as follows:
	\begin{align*}
		{ }_1 \Phi_1(a ; c ; x)=M(a ; c ; x)=\sum_{n=0}^{\infty} \frac{(a)_n}{(c)_n} \frac{x^n}{n !}
	\end{align*}
	where $|x|<\infty$; $c \neq 0,-1,-2, \ldots$. The confluent hypergeometric function of second kind is defined by 
	\begin{align*}
		U(a ; c ; x)=\frac{\pi}{\sin(\pi c)}\left[\frac{M(a ; c ; x)}{\Gamma(1+a-c) \Gamma(c)}-\frac{x^{1-c} M(1+a-c ; 2-c ; x)}{\Gamma(a) \Gamma(2-c)}\right].
	\end{align*}
	The integral representation of $U(a ; c ; x)$ is given by
	\begin{align*}
		U(a ; c ; x)=\frac{1}{\Gamma(a)} \int_0^{\infty} e^{-x t} t^{a-1}(1+t)^{c-a-1} d t
	\end{align*}
	where $ a>0, c>0, c \neq 1,2, \ldots$. In \cite{FS}, Ferreira and Salinas defined the incomplete generalized gamma function by using the confluent hypergeometric function of the second kind as follows:
	$$
	\text{ }_\lambda\gamma_\omega\left(p,\delta;a;c;\nu\right)=\int_0^\omega x^{\lambda-1}e^{-px^\delta}U(a;c;\nu x^\delta)dx
	$$
	where $x>0$, $\delta>0$, $p>0$, $a$ and $c$ are arbitrary constants. Motivated by this definition, we define the following integral transform
	\begin{align*}
		\text{ }_{\lambda}\gamma_\infty\left(p,\delta;a;c;\nu;f(x)\right)=\int_0^\infty x^{\lambda-1}e^{-px^\delta}U(a;c;\nu x^\delta)f(x)dx.
	\end{align*}
	In \cite{FS}, Ferreira and Salinas evaluated the following integral,
	\begin{align}
		\int_0^{\infty} x^{\lambda-1} e^{-p x^\delta} U\left(a ; c ; v x^\delta\right) d x&=\frac{\pi}{\delta \sin(\pi c) p^{\frac{\lambda}{\delta}}}{\left[\frac{A(p, \delta, \lambda, a, c, v)-B(p, \delta, \lambda, a, c, v)}{p^{1-c} \Gamma(1+a-c) \Gamma(c) \Gamma(a) \Gamma(2-c)}\right]}\label{UI}
	\end{align}
	where
	\begin{align*}
		& A(p, \delta, \lambda, a, c, v)=p^{1-c} \Gamma(a) \Gamma\left(\frac{\lambda}{\delta}\right) \Gamma(2-c) { }_2 F_1\left(a, \frac{\lambda}{\delta} ; c ; \frac{v}{p}\right),\\
		& B(p, \delta, \lambda, a, c, v)=v^{1-c} \Gamma(1+a-c) \Gamma(c) \Gamma\left(\frac{\lambda}{\delta}-c+1\right)
		{ }_2 F_1\left(1+a-c, \frac{\lambda}{\delta}-c+1 ; 2-c ; \frac{v}{p}\right)
	\end{align*}
	where $\lambda, \delta, p>0$, $\lambda, v, p$ are constants such as $0<v<p, c<1$, $c \notin \mathbb{Z}$ and $a, 1+a-c \notin \mathbb{Z}^{-}$. But, using the relation
	\begin{align*}
		{ }_2 F_1(a, b ; c ; z)&=\frac{\Gamma(a-c+1) \Gamma(b-c+1)}{\Gamma(1-c) \Gamma(a+b-c+1)}{ }_2 F_1(a, b ; a+b-c+1 ; 1-z)\\
		&-
		\frac{\Gamma(a-c+1) \Gamma(b-c+1) \Gamma(c-1)}{\Gamma(a) \Gamma(b) \Gamma(1-c)} z^{1-c}{ }_2 F_1(a-c+1, b-c+1 ; 2-c ; z)
	\end{align*}
	of the generalized hypergeometric function, we can write it as
	\begin{align}
		\int_0^{\infty} x^{\lambda-1} e^{-p x^\delta} U\left(a ; c ; v x^\delta\right) d x&=\frac{1}{\delta p^{\frac{\lambda}{\delta}}}\frac{\Gamma\left(\frac{\lambda}{\delta}\right)\Gamma\left(\frac{\lambda}{\delta}-c+1\right)}{\Gamma\left(a+\frac{\lambda}{\delta}-c+1\right)}{ }_2 F_1\left(a, \frac{\lambda}{\delta};a+\frac{\lambda}{\delta}-c+1 ;1-\frac{v}{p}\right).\label{UII}
	\end{align}
	
	\section{Parseval-Goldstein Type Theorems}
	
	In this section, we will prove some identities and Parseval-Goldstein type theorems. 
	
	The following lemma shows that the generalized Stieltjes transform can be obtained by applying $\mathcal{L}_{\alpha,\mu}$-transform and $\mathcal{L}_{\delta,\mu}$-transform consecutively.
	
	\begin{lemma}
		\label{Lmm1} Let $F\left(  y\right)  =\mathcal{L}_{\delta,\mu}\left\{  f\left(  t\right)
		;y\right\}$. If $x,y,\alpha,\delta\in\mathbb{C}$, $\mu\in\mathbb{R}$ and $f,F\in A$, then the following identity%
		\begin{equation}
			\mathcal{L}_{\alpha,\mu}\left\{\mathcal{L}_{\delta,\mu}\left\{  f\left(  t\right)  ;x\right\}  ;y\right\}
			=\dfrac{1}{\mu}\Gamma\left(  \frac{\alpha}{\mu
			}\right)  \mathcal{S}_{\delta,\mu,\frac{\alpha}{\mu}}\left\{  f\left(  t\right)  ;y\right\}  \label{LLS}%
		\end{equation}
		holds true for $\operatorname{Re}\alpha> \mu>0$, $\operatorname{Re}\delta> \mu>0$, $\operatorname{Re}y>0$, $\operatorname{Re}x>0$, $\operatorname{Re}\left(\frac{\alpha}{\mu}\right)>0$ provided that the integrals involved converge absolutely.
	\end{lemma}
	
	\begin{proof}
		Using the definition of (\ref{LAM}), changing the order of integration, which
		is permissible by absolute convergence of the integrals involved, we get
		\begin{align*}
			\mathcal{L}_{\alpha,\mu}\left\{\mathcal{L}_{\delta,\mu}\left\{  f\left(  t\right)  ;x\right\}  ;y\right\}=\int_0^\infty t^{\delta-1}f(t)\mathcal{L}_{\alpha,\mu}\{1;\sqrt[\mu]{t^\mu+y^\mu}\}dt
		\end{align*}
		and using the relation (\ref{LAML}) and the formula
		\begin{align*}
			\mathcal{L}_{\alpha,\mu}\{1;\sqrt[\mu]{t^\mu+y^\mu}\}=\Gamma\left(\frac{\alpha}{\mu}\right)\frac{1}{\mu}\frac{1}{(t^\mu+y^\mu)^{\frac{\alpha}{\mu}}},
		\end{align*}
		we arrive at (\ref{LLS}).
	\end{proof}
	
	The following is a Parseval-Goldstein type theorem for $\mathcal{L}_{\alpha,\mu}$-transform and generalized Stieltjes transform.
	\begin{theorem}
		\label{Thm1} If $f,g\in A$, $\alpha,\delta\in\mathbb{C}$, $\mu,y\in\mathbb{R}$ and  then the following identities%
		\begin{align}
			&  \int_{0}^{\infty}y^{\lambda-1}\mathcal{L}%
			_{\alpha,\mu}\left\{  f\left(  t\right)  ;y\right\}  \mathcal{L}_{\delta,\mu
			}\left\{  g\left(  x\right)  ;y\right\}  dy=\dfrac{1}{\mu}\Gamma\left(  \frac{\lambda}{\mu}\right)  \int_{0}^{\infty}t^{\alpha-1}f\left(  t\right)  \mathcal{S}%
			_{\delta,\mu,\frac{\lambda}{\mu}}\left\{  g\left(  x\right)  ;t\right\}  dt,\label{LLSAMR1}\\
			&  \int_{0}^{\infty}y^{\lambda-1}\mathcal{L}%
			_{\alpha,\mu}\left\{  f\left(  t\right)  ;y\right\}  \mathcal{L}_{\delta,\mu
			}\left\{  g\left(  x\right)  ;y\right\}  dy=\dfrac{1}{\mu}\Gamma\left(  \frac{\lambda
			}{\mu}\right)  \int_{0}^{\infty}x^{\delta-1}g\left(  x\right)  \mathcal{S}%
			_{\alpha,\mu,\frac{\lambda}{\mu}}\left\{  f\left(  t\right)  ;x\right\}  dx, \label{LLSAMR2}%
		\end{align}
		hold true for $\operatorname{Re}\alpha> \mu>0$, $\operatorname{Re}\delta> \mu>0$, $y>0$, $\operatorname{Re}\left(\frac{\lambda}{\mu}\right)>0$ provided that the integrals involved converge absolutely.
	\end{theorem}
	
	\begin{proof}
		Using the definition (\ref{LAM}) and changing the order of integration, we
		have%
		\[
		\int_{0}^{\infty}y^{\lambda-1}\mathcal{L}%
		_{\alpha,\mu}\left\{f\left(t\right);y\right\}  \mathcal{L}_{\delta,\mu
		}\left\{  g\left(  x\right)  ;y\right\}  dy=\int_{0}^{\infty}t^{\alpha
			-1}f\left(  t\right)  \mathcal{L}_{\lambda,\mu}\left\{\mathcal{L}_{\delta,\mu}\left\{  g\left(  x\right)
		;y\right\}  ;t\right\}  dt.
		\]
		Using the identity (\ref{LLS}) of Lemma \ref{Lmm1}, we arrive at
		(\ref{LLSAMR1}). Proof of (\ref{LLSAMR2}) is similar.
	\end{proof}
	
	As a result of Theorem \ref{Thm1}, the following relation can be obtained from the equivalence of relations (\ref{LLSAMR1}) and (\ref{LLSAMR2}).
	\begin{align}
		\int_{0}^{\infty}t^{\alpha-1}f\left(  t\right)  \mathcal{S}%
		_{\delta,\mu,\frac{\lambda}{\mu}}\left\{  g\left(  x\right)  ;t\right\}  dt= \int_{0}^{\infty}x^{\delta-1}g\left(  x\right)  \mathcal{S}%
		_{\alpha,\mu,\frac{\lambda}{\mu}}\left\{  f\left(  t\right)  ;x\right\}  dx.\label{SS1}%
	\end{align}
	
	The following lemma shows that the generalized Stieltjes transform can be obtained by applying $\mathcal{L}_{\alpha,\mu}$-transform and $\mathcal{F}_{s,\delta,\mu}$-transform consecutively or $\mathcal{L}_{\alpha,\mu}$-transform and $\mathcal{F}_{c,\delta,\mu}$-transform consecutively in both order.
	\begin{lemma}
		\label{Lmm2}  Let $F\left(  x\right)  =\mathcal{L}_{\alpha,\mu}\left\{
		f\left(  t\right)  ;x\right\}$, $F_s\left(  x\right)  =\mathcal{F}_{s,\delta,\mu}\left\{
		f\left(  t\right)  ;x\right\}$ and $F_c\left(  x\right)  =\mathcal{F}_{c,\delta,\mu}\left\{
		f\left(  t\right)  ;x\right\}$. If $f,F_s,F_c\in A$, $\alpha,\delta\in\mathbb{C}$, and $x,y,\mu\in\mathbb{R}$, then the following identities%
		\begin{align}
			\mathcal{L}_{\alpha,\mu}\left\{\mathcal{F}_{s,\delta,\mu}\left\{
			f\left(  t\right)  ;x\right\}  ;y\right\}&  =\dfrac{1}{\mu}\Gamma\left(\frac{\alpha}{\mu}\right)\mathcal{S}%
			_{\delta,2\mu,\frac{\alpha}{2\mu}}\left\{\sin\left[\frac{\alpha}{\mu}\arctan\left(\frac{t^\mu}{y^\mu}\right)\right] f\left(
			t\right);y\right\},\label{LAMFSAM}\\
			\mathcal{F}_{s,\delta,\mu}\left\{\mathcal{L}_{\alpha,\mu}\left\{
			f\left(  t\right)  ;x\right\}  ;y\right\}
			&=\dfrac{1}{\mu}\Gamma\left(\frac{\delta}{\mu}\right)\mathcal{S}%
			_{\alpha,2\mu,\frac{\delta}{2\mu}}\left\{\sin\left[\frac{\delta}{\mu}\arctan\left(\frac{y^\mu}{t^\mu}\right)\right] f\left(
			t\right);y\right\},\label{FSAMLAM}\\
			\mathcal{L}_{\alpha,\mu}\left\{\mathcal{F}_{c,\delta,\mu}\left\{
			f\left(  t\right)  ;x\right\}  ;y\right\}&  =\dfrac{1}{\mu}\Gamma\left(\frac{\alpha}{\mu}\right)\mathcal{S}%
			_{\delta,2\mu,\frac{\alpha}{2\mu}}\left\{\cos\left[\frac{\alpha}{\mu}\arctan\left(\frac{t^\mu}{y^\mu}\right)\right] f\left(
			t\right);y\right\},\label{LAMFCAM}\\
			\mathcal{F}_{c,\delta,\mu}\left\{\mathcal{L}_{\alpha,\mu}\left\{
			f\left(  t\right)  ;x\right\}  ;y\right\}
			&=\dfrac{1}{\mu}\Gamma\left(\frac{\delta}{\mu}\right)\mathcal{S}%
			_{\alpha,2\mu,\frac{\delta}{2\mu}}\left\{\cos\left[\frac{\delta}{\mu}\arctan\left(\frac{y^\mu}{t^\mu}\right)\right] f\left(
			t\right);y\right\}, \label{FCAMLAM}
		\end{align}
		hold true for  $\operatorname{Re}\alpha> \mu>0$, $x^\mu>0$, $y^\mu>0$, $\operatorname{Re}\left(\frac{\alpha}{2\mu}\right)>0$, $\operatorname{Re}\left(\frac{\delta}{2\mu}\right)>0$, $f(t)$ and $F(x)$ are piecewise continuous, $t^{\delta-1}f(t)$ and $x^{\delta-1} F(x)$ are absolutely integrable over $[0,\infty)$ provided that the integrals involved converge absolutely.
	\end{lemma}
	
	\begin{proof}
		Using the definitions of (\ref{LAM}) and (\ref{FSAM}), changing the order of
		integration, which is permissible by absolute convergence of the integrals
		involved, we have%
		\[
		\mathcal{L}_{\alpha,\mu}\left\{\mathcal{F}_{s,\delta,\mu}\left\{  f\left(
		t\right)  ;x\right\}  ;y\right\}  =\int_{0}^{\infty}t^{\delta-1}f\left(
		t\right)  \mathcal{L}_{\alpha,\mu}\left\{  \sin\left(  x^{\mu}t^{\mu}\right)
		;y\right\}  dt.
		\]
		Using the relation (\ref{LAML}), the known formula \cite[p.152, Entry(15)]
		{E1}
		\[
		\mathcal{L}\left\{t^{\nu-1} \sin\left( at\right)  ;y\right\}=\dfrac{\Gamma\left(a\right)}{\left(a^2+y^2\right)^{\nu/2}}\sin\left[\nu\arctan\left(\frac{a}{y}\right)\right]
		\]
		where $\operatorname{Re}(\nu)>-1$, $\operatorname{Re}(y)>|\operatorname{Im}(a)|$ and definition of (\ref{SAMR}), we arrive at (\ref{LAMFSAM}). Similarly, using
		the definitions of (\ref{LAM}) and (\ref{FCAM}) changing the order of
		integration, which is permissible by absolute convergence of the integrals
		involved, we have%
		\[
		\mathcal{L}_{\alpha,\mu}\left\{\mathcal{F}_{c,\delta,\mu}\left\{  f\left(
		t\right)  ;x\right\}  ;y\right\}=\int_{0}^{\infty}t^{\delta-1}f\left(
		t\right)  \mathcal{L}_{\alpha,\mu}\left\{\cos\left(x^{\mu}t^{\mu}\right)
		;y\right\}  dt.
		\]
		Using the relation (\ref{LAML}), the known formula \cite[p.157, Entry(58)]{E1}
		\[
		\mathcal{L}\left\{t^{\nu-1} \cos\left( at\right)  ;y\right\}=\dfrac{\Gamma\left(a\right)}{\left(a^2+y^2\right)^{\nu/2}}\cos\left[\nu\arctan\left(\frac{a}{y}\right)\right]
		\]
		where $\operatorname{Re}(\nu)>0$, $\operatorname{Re}(y)>|\operatorname{Im}(a)|$ and definition of (\ref{SAMR}), we arrive at (\ref{LAMFCAM}). Proof of
		(\ref{FSAMLAM}) and (\ref{FCAMLAM}) are similar and can be made using the same
		definitions, relations and formulas.
	\end{proof}
	
	The following is a Parseval-Goldstein type theorem for $\mathcal{L}_{\alpha,\mu}$-transform, generalized Fourier cosine and sine transforms and generalized Stieltjes transform.
	
	\begin{theorem}
		\label{Thm2} If $f\in A$, $g(x)$ is piecewise continuous and $t^{\delta-1}f(t)$ is absolutely integrable over $[0,\infty)$, $\alpha\in\mathbb{C}$, $\mu,y\in\mathbb{R}$, then the following identities%
		\begin{align}
			&  \int_{0}^{\infty}y^{\lambda-1}\mathcal{L}_{\alpha,\mu}\left\{  f\left(
			t\right)  ;y\right\}  \mathcal{F}_{s,\delta,\mu}\left\{  g\left(  x\right)
			;y\right\}dy\nonumber\\
			&\quad=\dfrac{1}{\mu}\Gamma\left(\frac{\lambda}{\mu}\right)\int_{0}^{\infty
			}t^{\alpha-1}f\left(  t\right)  \mathcal{S}%
			_{\delta,2\mu,\frac{\lambda}{2\mu}}\left\{\sin\left[\frac{\lambda}{\mu}\arctan\left(\frac{x^\mu}{t^\mu}\right)\right] g\left(
			x\right);t\right\}
			dt,\label{LFSS1}\\
			&\int_{0}^{\infty}y^{\lambda-1}\mathcal{L}_{\alpha,\mu}\left\{  f\left(
			t\right)  ;y\right\}  \mathcal{F}_{s,\delta,\mu}\left\{  g\left(  x\right)
			;y\right\}dy\nonumber\\
			&\quad=\dfrac{1}{\mu}\Gamma\left(\frac{\lambda}{\mu}\right)\int_{0}^{\infty
			}x^{\delta-1}g\left(  x\right)  \mathcal{S}%
			_{\delta,2\mu,\frac{\lambda}{2\mu}}\left\{\sin\left[\frac{\lambda}{\mu}\arctan\left(\frac{x^\mu}{t^\mu}\right)\right] f\left(
			t\right);x\right\}
			dx,\label{LFSS2}\\
			&\int_{0}^{\infty}y^{\lambda-1}\mathcal{L}_{\alpha,\mu}\left\{  f\left(
			t\right)  ;y\right\}  \mathcal{F}_{c,\delta,\mu}\left\{  g\left(  x\right)
			;y\right\}  dy\nonumber\\
			&\quad=\dfrac{1}{\mu}\Gamma\left(\frac{\lambda}{\mu}\right)\int_{0}^{\infty
			}t^{\alpha-1}f\left(  t\right)  \mathcal{S}%
			_{\delta,2\mu,\frac{\lambda}{2\mu}}\left\{\cos\left[\frac{\lambda}{\mu}\arctan\left(\frac{x^\mu}{t^\mu}\right)\right] g\left(
			x\right);t\right\}
			dt,\label{LFCS1}\\
			&\int_{0}^{\infty}y^{\lambda-1}\mathcal{L}_{\alpha,\mu}\left\{  f\left(
			t\right)  ;y\right\}  \mathcal{F}_{c,\delta,\mu}\left\{  g\left(  x\right)
			;y\right\}  dy\nonumber\\
			&\quad=\dfrac{1}{\mu}\Gamma\left(\frac{\lambda}{\mu}\right)\int_{0}^{\infty
			}x^{\delta-1}g\left(  x\right)  \mathcal{S}%
			_{\delta,2\mu,\frac{\lambda}{2\mu}}\left\{\cos\left[\frac{\lambda}{\mu}\arctan\left(\frac{x^\mu}{t^\mu}\right)\right] f\left(
			t\right);x\right\}
			dx,
			\label{LFCS2}%
		\end{align}
		hold true for $\operatorname{Re}\alpha> \mu>0$, $y^\mu>0$, $\operatorname{Re}\left(\frac{\lambda}{2\mu}\right)>0$ provided that the integrals involved converge absolutely.
	\end{theorem}
	
	\begin{proof}
		Using the definition (\ref{LAM}) and changing the order of integration, we
		have%
		\begin{align*}
			&\int_{0}^{\infty}y^{\lambda-1}\mathcal{L}_{\alpha,\mu}\left\{  f\left(
			t\right)  ;y\right\}  \mathcal{F}_{s,\delta,\mu}\left\{  g\left(  x\right)
			;y\right\}  dy=\int_{0}^{\infty}t^{\alpha-1}f\left(  t\right)  \mathcal{L}%
			_{\lambda,\mu}\left\{\mathcal{F}_{s,\delta,\mu}\left\{  g\left(  x\right)
			;y\right\}  ;t\right\}  dt.
		\end{align*}
		Using the identity (\ref{LAMFSAM}) of Lemma (\ref{Lmm2}), we arrive at
		(\ref{LFSS1}). Proof of (\ref{LFSS2}) is similar and can be made using the
		definition (\ref{FSAM}) and identity (\ref{FSAMLAM}) of Lemma (\ref{Lmm2}).
		Using the definition (\ref{LAM}), changing the order of integration and using
		the identity (\ref{LAMFCAM}) of Lemma (\ref{Lmm2}), we arrive at
		(\ref{LFCS1}). Proof of (\ref{LFCS2}) is similar and can be made using the
		definition (\ref{FCAM}) and identity (\ref{FCAMLAM}) of Lemma (\ref{Lmm2}).
	\end{proof}
	
	The following lemma shows that the improper integral involving the confluent hypergeometric function of second kind can be obtained by applying $\mathcal{L}_{\alpha,\mu}$-transform and generalized Stieltjes integral transform consecutively in both order.
	
	\begin{lemma}
		\label{Lmm3} Let $S\left(  x\right)  =\mathcal{S}_{\delta,\mu, \rho}\left\{  f\left(  t\right)
		;x\right\}$. If $\alpha\in\mathbb{C}$, $x,y,\mu\in\mathbb{R}$ and $f,S\in A$, then the following identities
		\begin{align}
			\mathcal{L}_{\alpha,\mu}\left\{	\mathcal{S}_{\delta,\mu,\rho}\left\{f(t);x\right\};y\right\}=\frac{\Gamma(\frac{\alpha}{\mu})}{\mu}\int_0^\infty t^{\delta+\alpha-\mu\rho-1}U\left(\frac{\alpha}{\mu};1+\frac{\alpha}{\mu}-\rho;t^\mu y^\mu\right)f(t)dt,\label{LAMSDMR}\\
			\mathcal{S}_{\delta,\mu,\rho}\left\{\mathcal{L}_{\alpha,\mu}\left\{f(t);x\right\};y\right\}=\frac{y^{\delta-\mu\rho}}{\mu}\Gamma\left(\frac{\delta}{\mu}\right)\int_0^\infty t^{\alpha-1}U\left(\frac{\delta}{\mu};1+\frac{\delta}{\mu}-\rho;t^\mu y^\mu\right)f(t)dt\label{SDMRLAM}
		\end{align}
		hold true for $\operatorname{Re}\alpha> \mu>0$, $\operatorname{Re}\delta>0$, $x>0$, $y>0$, $\operatorname{Re}\left(\frac{\alpha}{\mu}\right)>0$, $\operatorname{Re}\left(\frac{\delta}{\mu}\right)>0$, $\operatorname{Re}\left(1+\frac{\alpha}{\mu}\right)>\operatorname{Re}\rho>$, $\operatorname{Re}\left(1+\frac{\delta}{\mu}\right)>\operatorname{Re}\rho>0$ provided that the integrals involved converge absolutely.
	\end{lemma}
	\begin{proof}
		Using the definitions of (\ref{LAM}) and (\ref{SAMR}), changing the order of
		integration, which is permissible by absolute convergence of the integrals
		involved, we have%
		\[
		\mathcal{L}_{\alpha,\mu}\left\{\mathcal{S}_{\delta,\mu,\rho}\left\{  f\left(
		t\right)  ;x\right\}  ;y\right\} =\int_{0}^{\infty}t^{\delta-1}f\left(
		t\right) \left[\int_0^\infty \frac{x^{\alpha-1}e^{-x^\mu y^\mu}}{\left(t^\mu+x^\mu\right)^\rho}dx\right]dt.
		\]
		Now, making the change of variable $x=tu^{\frac{1}{\mu}}$ in the inner integral, we get
		\[
		\mathcal{L}_{\alpha,\mu}\left\{\mathcal{S}_{\delta,\mu,\rho}\left\{  f\left(
		t\right)  ;x\right\}  ;y\right\}=\frac{1}{\mu}\int_0^\infty t^{\delta+\alpha-\mu\rho-1}f(t)\left[\int_0^\infty\frac{u^{\frac{\alpha}{\mu}-1}}{(1+u)^\rho}e^{-t^\mu y^\mu u}du\right]dt.
		\]
		Using the integral representation of the confluent hypergeometric function $U(a, b, z)$, we arrive at (\ref{LAMSDMR}). Proof of (\ref{SDMRLAM}) is similar and can be made using the same definitions and formulas.
	\end{proof}
	
	The following is a Parseval-Goldstein type theorem for $\mathcal{L}_{\alpha,\mu}$-transform,
	generalized Stieltjes transform and $\text{ }_{\lambda}\gamma_\infty-$transform.
	\begin{theorem}
		\label{Thm3} If $\alpha,\delta,\rho\in\mathbb{C}$, $\mu\in\mathbb{R}$ and $f\in A$, then the following identities
		\begin{align}
			&  \int_{0}^{\infty}y^{\lambda-1}\mathcal{L}_{\alpha,\mu}\left\{  f\left(
			t\right)  ;y\right\}  \mathcal{S}_{\delta,\mu,\rho}\left\{  g\left(  x\right)
			;y\right\}dy\nonumber\\
			&\quad=\frac{\Gamma(\frac{\lambda}{\mu})}{\mu}\int_{0}^{\infty
			}t^{\alpha-1}f\left( t\right) \text{ }_{\delta+\lambda-\mu\rho}\gamma_\infty\left(0;\mu;\frac{\lambda}{\mu};1+\frac{\lambda}{\mu}-\rho;t^\mu;g(x)\right)dt,\label{LFSR1}\\
			&\int_{0}^{\infty}y^{\lambda-1}\mathcal{L}_{\alpha,\mu}\left\{  f\left(
			t\right)  ;y\right\}  \mathcal{S}_{\delta,\mu,\rho}\left\{  g\left(  x\right)
			;y\right\}dy\nonumber\\
			&\quad=\frac{\Gamma(\frac{\lambda}{\mu})}{\mu}\int_{0}^{\infty
			}x^{\delta+\lambda-\mu\rho-1}g\left(  x\right) \text{ }_{\alpha}\gamma_\infty\left(0;\mu;\frac{\lambda}{\mu};1+\frac{\lambda}{\mu}-\rho;x^\mu;f(t)\right)dx,\label{LFSR2}
		\end{align}
		hold true for $\operatorname{Re}\alpha> \mu>0$, $\operatorname{Re}\delta>0$, $y>0$, $\operatorname{Re}\left(\frac{\alpha}{\mu}\right)>0$, $\operatorname{Re}\left(\frac{\delta}{\mu}\right)>0$, $\operatorname{Re}\left(1+\frac{\alpha}{\mu}\right)>\operatorname{Re}\rho>$, $\operatorname{Re}\left(1+\frac{\delta}{\mu}\right)>\operatorname{Re}\rho>0$ provided that the integrals involved converge absolutely.
	\end{theorem}
	\begin{proof} Using the definition (\ref{LAM}) and changing the order of integration, we
		have%
		\begin{align*}
			&\int_{0}^{\infty}y^{\lambda-1}\mathcal{L}_{\alpha,\mu}\left\{  f\left(
			t\right)  ;y\right\}  \mathcal{S}_{\delta,\mu,\rho}\left\{  g\left(  x\right)
			;y\right\}dy=\int_{0}^{\infty}t^{\alpha-1}f\left(  t\right)  \mathcal{L}%
			_{\lambda,\mu}\left\{\mathcal{S}_{\delta,\mu,\rho}\left\{  g\left(  x\right)
			;y\right\}  ;t\right\}  dt.
		\end{align*}
		Using the identity (\ref{LAMSDMR}) of Lemma \ref{Lmm3}, we arrive at
		(\ref{LFSR1}). Proof of (\ref{LFSR2}) is similar and can be using the
		definition (\ref{SAMR}) and identity (\ref{SDMRLAM}) of Lemma \ref{Lmm3}.
	\end{proof}
	
	\section{Applications}
	We know that \cite{A1}
	\begin{align}
		\mathcal{L}_{\alpha,\mu}\left\{x^{\lambda-1};y\right\}=\dfrac{1}{\mu}\Gamma\left(\dfrac{\alpha+\lambda-1}{\mu}\right)\dfrac{1}{y^{\alpha+\lambda-1}}\label{app0}
	\end{align}
	where $\operatorname{Re}y>0$ and $\operatorname{Re}\left(\frac{\alpha+\lambda-1}{\mu}\right)>-1$. 
	
	In this section we give some applications of above lemmas and theorems.  
	\begin{example}
		We show that
		\begin{align}
			\mathcal{S}_{\delta,\mu,\frac{\alpha}{\mu}}\left\{t^{\lambda-1};y\right\}&=\frac{1}{\mu y^{\alpha-\delta-\lambda+1}}B\left(\dfrac{\delta+\lambda-1}{\mu},\dfrac{\alpha-\delta-\lambda+1}{\mu}\right),\label{app1}\\
			\mathcal{S}_{\delta,\mu,\frac{\alpha}{\mu}}\left\{e^{-a^\mu t^\mu};y\right\}&=\frac{a^{\alpha-\delta}}{\mu}\Gamma\left(\dfrac{\delta}{\mu}\right)U\left(\dfrac{\alpha}{\mu};1+\dfrac{\alpha}{\mu}-\dfrac{\delta}{\mu};a^\mu y^\mu\right)\label{app2}
		\end{align}
		where $\operatorname{Re}\left(\frac{\delta+\lambda-1}{\mu}\right)>0$, $\operatorname{Re}\left(\frac{\alpha-\delta-\lambda+1}{\mu}\right)>0$, $\operatorname{Re}\left(\frac{\delta}{\mu}\right)>0$, $\operatorname{Re}\left(\frac{\alpha}{\mu}\right)>0$, $\operatorname{Re}\left(\frac{\alpha}{\mu}-\frac{\delta}{\mu}\right)>-1$ and $U$ is a second kind of confluent hypergeometric function.
	\end{example} 
	Setting $f(t)=t^{\lambda-1}$ in (\ref{LLS}), we have
	\begin{align*}
		\mathcal{S}_{\delta,\mu,\frac{\alpha}{\mu}}\left\{  t^{\lambda-1}  ;y\right\}=\mu\left\{\Gamma\left(  \frac{\alpha}{\mu
		}\right)  \right\}^{-1}\mathcal{L}_{\alpha,\mu}\left\{\mathcal{L}_{\delta,\mu}\left\{  t^{\lambda-1}  ;x\right\}  ;y\right\}.
	\end{align*}
	Using the formula (\ref{app0}) successively, we obtain the formula (\ref{app1}). Setting $f(t)=e^{-a^\mu t^\mu}$ in (\ref{LLS}), we have
	\begin{align*}
		\mathcal{S}_{\delta,\mu,\frac{\alpha}{\mu}}\left\{ e^{-a^\mu t^\mu}  ;y\right\}&=\mu\left\{\Gamma\left(  \frac{\alpha}{\mu
		}\right)  \right\}^{-1}\mathcal{L}_{\alpha,\mu}\left\{\mathcal{L}_{\delta,\mu}\left\{  e^{-a^\mu t^\mu}  ;x\right\}  ;y\right\}.
	\end{align*}
	Using the definition (\ref{LAM}) and the formula (\ref{app0}) for $\lambda=1$, we have
	\begin{align*}
		\mathcal{S}_{\delta,\mu,\frac{\alpha}{\mu}}\left\{ e^{-a^\mu t^\mu}  ;y\right\}&=\mu\left\{\Gamma\left(  \frac{\alpha}{\mu
		}\right)  \right\}^{-1}\mathcal{L}_{\alpha,\mu}\left\{\mathcal{L}_{\delta,\mu}\left\{ 1 ;\sqrt[\mu]{a^\mu+x^\mu}\right\}  ;y\right\}\\
		&=\left\{\Gamma\left(  \frac{\alpha}{\mu
		}\right)  \right\}^{-1}\Gamma\left(  \frac{\delta}{\mu
		}\right)\mathcal{L}_{\alpha,\mu}\left\{\frac{1}{\left(a^\mu+x^\mu\right)^{\frac{\delta}{\mu}}} ;y\right\}\\
		&=\left\{\Gamma\left(  \frac{\alpha}{\mu
		}\right)  \right\}^{-1}\Gamma\left(  \frac{\delta}{\mu
		}\right)\int_{0}^{\infty}\frac{x^{\alpha-1}e^{-y^\mu x^\mu}}{\left(a^\mu+x^\mu\right)^{\frac{\delta}{\mu}}} dx.
	\end{align*}
	Now, making the change of variable $x=au^{\frac{1}{\mu}}$, we get
	\begin{align*}
		\mathcal{S}_{\delta,\mu,\frac{\alpha}{\mu}}\left\{ e^{-a^\mu t^\mu}  ;y\right\}
		&=\frac{a^{\alpha-\delta}}{\mu}\int_{0}^{\infty}u^{\frac{\alpha}{\mu}-1}\left(1+u\right)^{-\frac{\delta}{\mu}}e^{-y^\mu a^\mu u}dx.	
	\end{align*}
	Using the integral representation of the confluent hypergeometric function $U(a, b, z)$, we arrive at (\ref{app2}).
	
	\begin{example}
		We show that
		\begin{align}
			\int_{0}^{\infty}\frac{y^{\lambda-1}}{\left(a^\mu+y^\mu\right)^{\frac{\alpha}{\mu}}\left(b^\mu+y^\mu\right)^{\frac{\delta}{\mu}}}dy&=\frac{b^{\lambda-\delta}}{\mu a^\alpha}\frac{\Gamma\left(\frac{\lambda}{\mu}\right)\Gamma\left(\frac{\alpha-\lambda+\delta}{\mu}\right)}{\Gamma\left(\frac{\alpha+\delta}{\mu}\right)}{ }_2 F_1\left(\frac{\alpha}{\mu}, \frac{\lambda}{\mu};\frac{\alpha+\delta}{\mu} ;1-\frac{b^\mu}{a^\mu}\right).\label{app3}
		\end{align}
		where $\operatorname{Re}\left(\frac{\lambda}{\mu}\right)>0$, $\operatorname{Re}\left(\frac{\alpha}{\mu}\right)>0$, $\operatorname{Re}\left(\frac{\alpha+\lambda}{\mu}\right)>\operatorname{Re}\left(\frac{\lambda}{\mu}\right)>0$.
	\end{example}
	If we set $f(t)=e^{-a^\mu t^\mu}$ and $g(x)=e^{-b^\mu x^\mu}$ in (\ref{LLSAMR1}) and use the definition (\ref{LAM}) and the formula (\ref{app2}), we have
	\begin{align*}
		&\int_{0}^{\infty}y^{\lambda-1}\mathcal{L}_{\alpha,\mu}\left\{e^{-a^\mu t^\mu};y\right\}\mathcal{L}_{\delta,\mu}\left\{e^{-b^\mu x^\mu};y\right\}dy=\dfrac{1}{\mu}\Gamma\left(  \frac{\lambda}{\mu}\right)  \int_{0}^{\infty}t^{\alpha-1}e^{-a^\mu t^\mu}  \mathcal{S}%
		_{\delta,\mu,\frac{\lambda}{\mu}}\left\{  e^{-b^\mu x^\mu}  ;t\right\}  dt
	\end{align*}
	and
	\begin{align*}
		&\int_{0}^{\infty}y^{\lambda-1}\mathcal{L}_{\alpha,\mu}\left\{1;\sqrt[\mu]{a^\mu+y^\mu}\right\}\mathcal{L}_{\delta,\mu}\left\{1;\sqrt[\mu]{b^\mu+y^\mu}\right\}dy\\
		&\qquad=\dfrac{b^{\lambda-\delta}}{\mu^2}\Gamma\left(  \frac{\lambda}{\mu}\right)\Gamma\left(\dfrac{\delta}{\mu}\right)\int_{0}^{\infty}t^{\alpha-1}e^{-a^\mu t^\mu}  U\left(\frac{\lambda}{\mu};1+\frac{\lambda-\delta}{\mu};b^\mu t^\mu\right)  dt.
	\end{align*}
	Using the formulas (\ref{app0}) for $\lambda=1$ and (\ref{UII}), we arrive at (\ref{app3}).
	
	\begin{example}
		We show that
		\begin{align}
			&\mathcal{S}%
			_{\delta,2\mu,\frac{\alpha}{2\mu}}\left\{t^{-\nu}\sin\left[\frac{\alpha}{\mu}\arctan\left(\frac{t^\mu}{y^\mu}\right)\right] ;y\right\}=\frac{1}{y^{\alpha+\nu-\delta}\mu}B\left(\frac{\alpha-\delta+\nu}{\mu},\frac{\delta-\nu}{\mu}\right)\sin\left[\frac{\pi}{2}\left(\frac{\delta-\nu}{\mu}\right)\right]\label{app4}
		\end{align}
		where $0<\operatorname{Re}\left(\frac{\delta-\nu}{\mu}\right)<2$.
	\end{example}
	
	If we choose $f(t)=t^{-\nu}$ in (\ref{LAMFSAM}), we have
	\begin{align}
		\mathcal{L}_{\alpha,\mu}\left\{\mathcal{F}_{s,\delta,\mu}\left\{
		t^{-\nu}  ;x\right\}  ;y\right\}&  =\dfrac{1}{\mu}\Gamma\left(\frac{\alpha}{\mu}\right)\mathcal{S}%
		_{\delta,2\mu,\frac{\alpha}{2\mu}}\left\{\sin\left[\frac{\alpha}{\mu}\arctan\left(\frac{t^\mu}{y^\mu}\right)\right] t^{-\nu};y\right\}.\label{p331}
	\end{align} 
	Firstly, let's find the inner transform on the left side of the identity
	\begin{align*}
		\mathcal{F}_{s,\delta,\mu}\left\{
		t^{-\nu} ;x\right\}=\int_0^\infty t^{\delta-\nu-1}\sin\left(x^\mu t^\mu\right)dt.
	\end{align*}
	Making the change of variable $x=u^{\frac{1}{\mu}}$ and using the formula \cite[p.68, (1)]{E1}, we get
	\begin{align}
		\mathcal{F}_{s,\delta,\mu}\left\{
		t^{-\nu} ;x\right\}=\frac{x^ {\nu-\delta}}{\mu}\Gamma\left(\frac{\delta-\nu}{\mu}\right)\sin\left[\frac{\pi}{2}\left(\frac{\delta-\nu}{\mu}\right)\right].\label{p332}
	\end{align}
	Finally, setting the result (\ref{p332}) in (\ref{p331}) and using the formula (\ref{app0}), we arrive at (\ref{app4}).
	
	\begin{example}
		We show that
		\begin{align}
			&\mathcal{S}_{\lambda-\delta,2\mu,\frac{\alpha+\beta}{2\mu}}\left\{P^{-\nu}_{\frac{\alpha+\beta}{\mu}-1}\left[\frac{y^\mu}{\sqrt{y^{2\mu}+a^{2\mu}}}\right];a\right\}=\Gamma\left(\frac{\lambda-\delta}{\mu}\right)\frac{2^             {\frac{\alpha+\beta+\delta-\lambda}{\mu}-1}a^             {\lambda-\alpha-\beta-\delta-3\mu}\Gamma\left(\frac{\alpha+\beta+\delta-\lambda}{\mu}+\frac{\nu}{2}\right)}{\mu\Gamma\left(\frac{\alpha+\beta}{\mu}+\nu\right)\Gamma\left(\frac{\nu-1}{2}-\frac{\alpha+\beta+\delta-\lambda}{\mu}\right)\cos\left(\frac{\pi\delta}{2\mu}\right)}.\label{app5}
		\end{align}
		where $\operatorname{Re}\left(\frac{\alpha+\beta+\delta-\lambda}{\mu}+\frac{\nu}{2}\right)>0$, $\operatorname{Re}\left(\frac{\alpha+\beta}{\mu}+\nu\right)>0$ and  $\operatorname{Re}\left(\frac{\nu-1}{2}-\frac{\alpha+\beta+\delta-\lambda}{\mu}\right)>0$.
	\end{example}
	Setting $f(t)=t^{\beta}\mathcal{J}_\nu\left(a^\mu t^\mu\right)$ and $g(x)=1$ in (\ref{LFSS1}), we obtain
	\begin{align}
		&  \int_{0}^{\infty}y^{\lambda-1}\mathcal{L}_{\alpha,\mu}\left\{  t^{\beta}\mathcal{J}_\nu\left(a^\mu t^\mu\right)  ;y\right\}  \mathcal{F}_{s,\delta,\mu}\left\{  1
		;y\right\}dy\nonumber\\
		&=\dfrac{1}{\mu}\Gamma\left(\frac{\lambda}{\mu}\right)\int_{0}^{\infty
		}t^{\alpha+\beta-1}\mathcal{J}_\nu\left(a^\mu t^\mu\right)  \mathcal{S}%
		_{\delta,2\mu,\frac{\lambda}{2\mu}}\left\{\sin\left[\frac{\lambda}{\mu}\arctan\left(\frac{x^\mu}{t^\mu}\right)\right] ;t\right\}
		dt.\label{p341}
	\end{align}
	To start with, let's find the first transform on the left side of the identity
	\begin{align*}
		\mathcal{L}_{\alpha,\mu}\left\{ t^{\beta}\mathcal{J}_\nu\left(a^\mu t^\mu\right)  ;y\right\}=\frac{1}{\mu}\int_{0}^\infty t^{\alpha+\beta-1}e^{t^\mu y^\mu}\mathcal{J}_\nu\left(a^\mu t^\mu\right)dt.
	\end{align*}
	Making the change of variable $x=u^{\frac{1}{\mu}}$ and use the formula \cite[p.29, (6)]{E2}, we get
	\begin{align}
		\mathcal{L}_{\alpha,\mu}\left\{ t^{\beta}\mathcal{J}_\nu\left(a^\mu t^\mu\right)  ;y\right\}=\frac{1}{\mu}\Gamma\left(\frac{\alpha+\beta}{\mu}+\nu\right)\frac{1}{\left(y^{2\mu}+a^{2\mu}\right)^{\frac{\alpha+\beta}{2\mu}}}P^{-\nu}_{\frac{\alpha+\beta}{\mu}-1}\left[\frac{y^\mu}{\sqrt{y^{2\mu}+a^{2\mu}}}\right].\label{p342}
	\end{align}
	Secondly, using the formula (\ref{p332}) for $v=0$, we get
	\begin{align}
		\mathcal{F}_{s,\delta,\mu}\left\{1;y\right\}=\frac{1}{\mu}y^{-\delta}\sin\left(\frac{\pi\delta}{\mu}\right)\Gamma\left(\frac{\delta}{\mu}\right),\label{p343}
	\end{align}
	and using the formula (\ref{app4}), for $v=0$,	we obtain
	\begin{align}
		\mathcal{S}%
		_{\delta,2\mu,\frac{\lambda}{2\mu}}\left\{\sin\left[\frac{\lambda}{\mu}\arctan\left(\frac{x^\mu}{t^\mu}\right)\right] ;t\right\}=\frac{1}{t^{\lambda-\delta}\mu}B\left(\frac{\lambda-\delta}{\mu},\frac{\delta}{\mu}\right)\sin\left(\frac{\pi\delta}{2\mu}\right).\label{p344}
	\end{align}
	Now, setting the results (\ref{p342}), (\ref{p343}) and (\ref{p344}) in (\ref{p341}), we obtain 
	\begin{align*}
		&\Gamma\left(\frac{\alpha+\beta}{\mu}+\nu\right)\Gamma\left(\frac{\delta}{\mu}\right)\cos\left(\frac{\pi\delta}{2\mu}\right)\mathcal{S}_{\lambda-\delta,2\mu,\frac{\alpha+\beta}{2\mu}}\left\{P^{-\nu}_{\frac{\alpha+\beta}{\mu}-1}\left[\frac{y^\mu}{\sqrt{y^{2\mu}+a^{2\mu}}}\right];a\right\}\\
		&=\dfrac{1}{2}\Gamma\left(\frac{\lambda}{\mu}\right)B\left(\frac{\lambda-\delta}{\mu},\frac{\delta}{\mu}\right)\int_{0}^{\infty
		}t^{\alpha+\beta+\delta-\lambda-1}\mathcal{J}_\nu\left(a^\mu t^\mu\right)dt.
	\end{align*}
	Finally, making the change of variable $x=u^{\frac{1}{\mu}}$ on the right side of the equation and using the formula \cite[p.22, (7)]{E2}, we arrive at (\ref{app5}).
	
	\begin{example}
		We show that
		\begin{align}
			&\mathcal{S}_{\lambda+\delta-\mu\rho,\mu,\frac{\alpha}{\mu}}\left\{U\left(\frac{\delta}{\mu};1+\frac{\delta}{\mu}-\rho;b^\mu y^\mu\right);a\right\}=\frac{1}{a^{\alpha+\lambda-\mu\rho}}\frac{\Gamma\left(\rho-\frac{\lambda}{\mu}\right)}{\Gamma\left(\rho\right)\Gamma\left(\frac{\alpha}{\mu}\right)\Gamma\left(\frac{\delta}{\mu}\right)\Gamma\left(1+\rho-\frac{\lambda}{\mu}\right)}\nonumber\\
			&\quad\times\left\{\frac{a^{\mu\rho-\lambda}}{b^\delta\mu}\Gamma\left(\frac{\lambda}{\mu}\right)\Gamma\left(\frac{\delta}{\mu}\right)\Gamma\left(\frac{\alpha}{\mu}+\frac{\lambda}{\mu}-\rho\right)\Gamma\left(1+\rho-\frac{\lambda}{\mu}\right)\text{ }_3F_1\left(\frac{\delta}{\mu},\frac{\lambda}{\mu},\frac{\delta}{\mu}+\frac{\lambda}{\mu}-\rho;1+\frac{\lambda}{\mu}-\rho;\frac{b^\mu}{a^\mu}\right)\right.\nonumber\\
			&\quad\left.+\frac{1}{b^{\rho-\frac{\lambda}{\mu}+\frac{\alpha}{\mu}}\mu}\Gamma\left(\rho\right)\Gamma\left(\frac{\alpha}{\mu}\right)\Gamma\left(\frac{\lambda}{\mu}+1-\rho\right)\Gamma\left(\frac{\alpha}{\mu}-\frac{\lambda}{\mu}+\rho\right)\text{ }_3F_1\left(\frac{\alpha}{\mu}-\frac{\lambda}{\mu}+\rho,\rho,\frac{\alpha}{\mu}
			;1+\rho-\frac{\lambda}{\mu};\frac{b^\mu}{a^\mu}\right)\right\}.\label{app6}
		\end{align}
		where $\operatorname{Re}\frac{\alpha}{\mu}>0$, $\operatorname{Re}\frac{\delta}{\mu}>0$, $\operatorname{Re}\frac{\lambda}{\mu}>0$ and $\operatorname{Re}\left(\frac{\alpha+\lambda}{\mu}\right)>\operatorname{Re}\rho>\operatorname{Re}\frac{\lambda}{\mu}-1$.
	\end{example}
	Setting $f(t)=e^{-a^\mu t^\mu}$ and $g(x)=e^{-b^\mu x^\mu}$ in (\ref{LFSR2}), we obtain
	\begin{align*}
		&  \int_{0}^{\infty}y^{\lambda-1}\mathcal{L}_{\alpha,\mu}\left\{  e^{-a^\mu t^\mu}  ;y\right\}  \mathcal{S}_{\delta,\mu,\rho}\left\{ e^{-b^\mu t^\mu}
		;y\right\}dy\\
		&\qquad=\frac{\Gamma(\frac{\lambda}{\mu})}{\mu}\int_{0}^{\infty
		}x^{\delta-1}e^{-b^\mu x^\mu} \text{ }_{\alpha+\lambda-\mu\rho}\gamma_\infty\left(0;\mu;\frac{\lambda}{\mu};1+\frac{\lambda}{\mu}-\rho;x^\mu;e^{-a^\mu t^\mu}\right)dx.
	\end{align*}
	On the other hand, using the formulas (\ref{app2}), (\ref{UI}), (\ref{FI}) and
	\begin{align*}
		\mathcal{L}_{\alpha,\mu}\left\{  e^{-a^\mu t^\mu}  ;y\right\}=\mathcal{L}_{\alpha,\mu}\left\{ 1  ;\sqrt[\mu]{y^\mu+a^\mu}\right\}=\Gamma\left(\frac{\alpha}{\mu}\right)\frac{1}{\mu}\frac{1}{(a^\mu+y^\mu)^{\frac{\alpha}{\mu}}}
	\end{align*}
	and the definition (\ref{SAMR}), we arrive at (\ref{app6}).

	\section{Conclusion}
	
	In this work, we establish Parseval-Goldstein type relations and identities that include various integral transforms such as $\mathcal{L}_{\alpha,\mu}$-transform and generalized Stieltjes transform. Thus, using these results, we show how simple it can be to evaulate integral transforms of some elementary and special functions. It is possible to obtain all the results and applications in \cite{Y1,Y2,YS} when $\alpha=\mu=\delta=\lambda=\rho=1$ is chosen in all lemmas, theorems and applications in Sections 2 and 3.


\begin{thebibliography}{99}
		
	
	\bibitem{A1} Albayrak, D.: \emph{Theory and applications on a new generalized Laplace-type integral transform}, Math. Meth. Appl. Sci. \textbf{46}(4), 4363-4378 (2023).
	
	\bibitem{A2} Albayrak, D., Dernek, N.: \emph{Some Relations for the Generalized $ \widetilde{\mathcal{G}}_n, \widetilde{\mathcal{P}}_n$ Integral Transforms and Riemann-Liouville, Weyl Integral Operators}, Gazi Univ. J. Sci., \textbf{36}(1), 362 - 381 (2023).
	
	\bibitem{A3} Albayrak, D., Dernek, N.: \emph{On some generalized integral transforms and Parseval-Goldstein type relations}, Hacet. J. Math. Stat., \textbf{50}(2), 526-540 (2021).
	
	\bibitem{BD1} Al-Musallam, F., Kiryakova, V., Tuan,  V. K.: \emph{A multi-index Borel-Dzrbashjan transform}, Rocky Mt. J. Math., \textbf{32}(2), 409–428 (2002).
	
	\bibitem{BD2} Dzhrbashyan, M.M.: \emph{Integral transforms and representations of functions in the complex domain}, Moscow: Nauka, 1966.
	
	\bibitem{DB} Debnath, L., Bhatta, D.: \emph{Integral Transforms and Their Applications(3rd ed.)}, Chapman and Hall/CRC, 2014.
	
	\bibitem {E1} Erdélyi, A., Magnus, W.,  Oberhettinger, F., Tricomi, F. G.: \emph{Tables of Integral Transforms}, Vol. I, New York-Toronto-London: McGraw-Hill Book Company, Inc. 1954.
	
	\bibitem {E2} Erdélyi, A., Magnus, W.,  Oberhettinger, F., Tricomi, F. G.: \emph{Tables of Integral Transforms}, Vol. II, New York-Toronto-London: McGraw-Hill Book Company, Inc. 1954.
	
	\bibitem{FS} Ferreira, J., Salinas, S.: \emph{A gamma type distribution involving a confluent hypergeometric function of the second kind}, Rev. Téc. Ing. Univ. Zulia., \textbf{33}(2), 169-175  (2010).
	
	\bibitem{G} Glasser, M. L.: \emph{Some Bessel function integrals}, Kyungpook Math. J., \textbf{13}(2), 171–174 (1973).
	
	\bibitem{KKU} Karataş, H. B., Kumar, D., Uçar, F.: \emph{Some iteration and Parseval-Goldstein type identities with their applications}, Adv. Appl. Math. Sci., \textbf{29}(2), 563–574 (2020).
	
	\bibitem{KAU} Karataş, H. B., Albayrak, D., Uçar, F.: \emph{Some Parseval-Goldstein type identities with illustrative examples},  Proc. Inst. Math. Mech., \textbf{49}(1), 60-68 (2023).
	
	\bibitem{NVYO} Dernek, N., Kurt, V., Şimşek, Y.,  Yürekli, O.: \emph{A generalization of the Widder potential transform and applications, Integral Transforms and Special Functions}, \textbf{22}(6), 391-401 (2011).
	
	\bibitem{OSM} Oldham, K. B., Spanier, J., Myland,  J.: \emph{An atlas of functions}, Springer 2010.
	
	\bibitem{W} Widder, D. V.: \emph{A transform related to the Poisson integral for a half-plane}, Duke Math. J. \textbf{33}(2), 355–362 (1966).
	
	\bibitem{Y1} Yürekli, O.: \emph{A Parseval-type theorem applied to certain integral transforms}, IMA J. Appl. Math., \textbf{42}(3), 241-249 (1989).
	
	\bibitem{Y2} Yürekli, O.: \emph{A theorem on the generalized Stieltjes transform}, J. Math. Anal. Appl., \textbf{168}(1), 63-71 (1992).
	
	\bibitem{YS} Yürekli, O., Sadek, I.: \emph{A Parseval Goldstein type theorem on the Widder potential and its applications}, Int. J. Math. Math. Sci. \textbf{14}, Article ID 160375, 517-524 (1991).
		
	\end{thebibliography}
\end{document}